\documentclass[12pt]{article}
\usepackage{amssymb,graphicx,rawfonts,t1enc}
\input{prepictex}
\input{pictex}
\input{postpictex}
\newcommand{\D}{\mathbf{d}}
\newcommand{\ZERO}{\mathbf{0}}
\newcommand{\ONE}{\mathbf{1}}
\newcommand{\wt}{\widetilde}
\newcommand{\R}{\mathbb R}
\newcommand{\C}{\mathbb C}
\newcommand{\Q}{\mathbb Q}
\newcommand{\F}{\mathbb F}
\begin{document}
\hfill
\vskip 2.2truecm
\thispagestyle{empty}
\par
\noindent
{\large ~~~Mappings from ${\R}^n$ to ${\F}^n$ which preserve unit Euclidean distance,}
\par
\centerline{{\large where $\F$ is a field of characteristic $0$}}
\vskip 2.2truecm
\centerline{{\large Apoloniusz Tyszka}}
\vskip 2.2truecm
\par
\noindent
{\bf Summary.} Let $\F$ be a commutative field of characteristic $0$,
~$\varphi_n:{\F}^n \times {\F}^n \to \F$,
~$\varphi_n((x_1,...,x_n),$ $(y_1,...,y_n))=(x_1-y_1)^2+...+(x_n-y_n)^2$.
We say that $g:{\R}^n \to {\F}^n$ preserves distance $d \geq 0$ if
for each $x,y \in {\R}^n$ $|x-y|=d$ implies
$\varphi_n(g(x),g(y))=d^2$.
Let $f:{\R}^n \to {\F}^n$ preserve unit distance.
We prove: {\bf (1)} if $n \geq 2$,
$x,y \in {\R}^n$ and $x \neq y$, then $\varphi_n(f(x),f(y)) \neq 0$,
{\bf (2)} if $A,B,C,D \in {\R}^2$, $r \in \Q$ and
$\overrightarrow{CD}=r\overrightarrow{AB}$,
then $\overrightarrow{f(C)f(D)}=r\overrightarrow{f(A)f(B)}$,
{\bf (3)} if $A,B,C,D \in {\R}^2$ and $\overrightarrow{AB}$ and
$\overrightarrow{CD}$ are linearly dependent,
then $\overrightarrow{f(A)f(B)}$ and $\overrightarrow{f(C)f(D)}$
are linearly dependent,
{\bf (4)} if $A,B,C,D \in {\R}^2$ and
$\overrightarrow{AB}$ is perpendicular to $\overrightarrow{CD}$,
then $\overrightarrow{f(A)f(B)}$ is perpendicular to
$\overrightarrow{f(C)f(D)}$,
{\bf (5)} if $A,B,C,D \in {\R}^2$ and $|AB|=|CD|$,
then $\varphi_2(f(A),f(B))=\varphi_2(f(C),f(D))$.
Let $A_n(\F)$ denote the set of all
positive numbers $d$ such that any map $g:{\R}^n \to {\F}^n$
that preserves unit distance preserves also distance $d$.
Let $D_{n}(\F)$ denote the set of all positive numbers~$d$
with the property:
if $x,y \in {\R}^n$ and $|x-y|=d$ then there exists a finite set $S_{xy}$
with $\left
\{x,y \right\} \subseteq S_{xy} \subseteq {\R}^n$ such that any map
map $g:S_{xy}\rightarrow {\F}^n$ that preserves unit distance
preserves also the distance between $x$ and $y$.
Obviously, $\{1\} \subseteq D_n(\F) \subseteq A_n(\F)$.
We prove: {\bf (6)}
$A_n(\C) \subseteq \left\{d>0: d^2 \in \Q \right\}$,
{\bf (7)}~$\left\{d>0: d^2 \in \Q \right\} \subseteq D_2(\F)$.
\newpage
Let $\F$ be a commutative field of characteristic $0$,
$\varphi_n: {\F}^n \times {\F}^n \to \F$,
$\varphi_n((x_1,...,x_n)$,
\par
\noindent
$(y_1,...,y_n))=(x_1-y_1)^2+...+(x_n-y_n)^2$.
We say that $f:{\R}^n \to {\F}^n$ preserves distance $d \geq 0$ if
for each $x,y \in {\R}^n$ $|x-y|=d$ implies
$\varphi_n(f(x),f(y))=d^2$.
In this paper we study unit-distance preserving mappings from
${\R}^n$ to ${\F}^n$.
\vskip 0.2truecm
\par
By a {\sl complex isometry} of ${\C}^n$ we understand any map
$h:{\C}^n \rightarrow {\C}^n$ of the form
$$
h(z_1,z_2,...,z_n)=(z_1',z_2',...,z_n')
$$
where
$$
z_j'=a_{0j}+a_{1j}z_{1}+a_{2j}z_{2}+...+a_{nj}z_{n} \hspace{0.5truecm}(j=1,2,...,n),
$$
the coefficients $a_{ij}$ are complex and the matrix $||a_{ij}||$
$(i,j=1,2,...,n)$ is orthogonal i.e. satisfies the condition
$$
\sum_{j=1}^{n} a_{\mu j} a_{\nu j}=\delta_{\nu}^{\mu} \hspace{0.5truecm} (\mu,\nu=1,2,...,n)
$$
with Kronecker's delta. According to \cite{Borsuk},
$\varphi_n(x,y)$ is invariant under complex isometries i.e. for every
complex isometry $h:{\C}^n \rightarrow {\C}^n$
\vskip 0.2truecm
\par
\noindent
{\boldmath $\left(\diamond\right)$} \hspace{4.0truecm} $\forall x,y \in {\C}^n ~~\varphi_n(h(x),h(y))=\varphi_n(x,y).$
\vskip 0.2truecm
\par
\noindent
Conversely, if $h:{\C}^n \to {\C}^n$ satisfies
{\boldmath $\left(\diamond\right)$} then $h$ is a complex isometry;
it follows from [4, proposition 1, page 21] by replacing $\R$ with $\C$
and $d(x,y)$ with $\varphi_n(x,y)$. Similarly, if
$f:{\R}^n \to {\C}^n$ preserves all distances, then there exists
a complex isometry $h:{\C}^n \to {\C}^n$ such that
$f=h_{|{\R}^n}$.
\vskip 0.2truecm
\par
By a field endomorphism of $\C$ we understand any map
$g:\C \to \C$ satisfying:
\par
\centerline
{$\forall x,y \in {\C}~~g(x+y)=g(x)+g(y)$,}
\centerline{$\forall x,y \in {\C}~~g(x \cdot y)=g(x) \cdot g(y)$,}
\centerline{$g(0)=0$,}
\centerline{$g(1)=1$.}
\par
\noindent
Bijective endomorphisms are called automorphisms,
for more information on field endomorphisms and automorphisms
of $\C$ the reader is referred to \cite{Kestelman}, \cite{Kuczma} and \cite{Yale}.
If $r$ is a rational number, then $g(r)=r$ for any field
endomorphism $g:\C \to \C$.
Proposition~1 shows that only rational numbers $r$ have this property:
\vskip 0.2truecm
\par
\noindent
{\bf Proposition 1} (\cite{Tyszka200?}).
If $r \in \C$ and $r$ is not a rational number,
then there exists a field automorphism $g:\C \to \C$ such that
$g(r) \neq r$.
\vskip 0.2truecm
\par
\noindent
{\it Proof.} Note first that if $E$ is any subfield of $\C$
and if $g$ is an automorphism of $E$, then~$g$ can be extended
to an automorphism of $\C$. This follows from [7, corollary 1, A.V.111]
by taking $\Omega=\C$ and $K=\Q$.
Now let $r \in \C \setminus \Q$.
If $r$ is algebraic over $\Q$, let $E$ be the
splitting field in $\C$ of the minimal polynomial $\mu$ of $r$ over $\Q$
and let $r' \in E$ be any other root of $\mu$.
Then there exists an automorphism $g$ of $E$ that sends $r$ to $r'$,
see for example [9, corollary 2, page 66].
If $r$ is transcendental over $\Q$, let $E={\Q}(r)$ and let
$r' \in E$ be any other generator of $E$ (e.g. $r'=1/r$). Then there
exists an automorphism~$g$ of $E$ that sends $r$ to $r'$.
In each case, $g$ can be extended to an automorphism of $\C$.
\vskip 0.2truecm
\par
If $g:\C \to \C$ is a field endomorphism then
$(g_{|\R},...,g_{|\R}): {\R}^n \to {\C}^n$ preserves all distances
$\sqrt{r}$ with rational $r \geq 0$.
Indeed, if $(x_1-y_1)^2+...+(x_n-y_n)^2=(\sqrt{r})^2$ then
\begin{eqnarray*}
\varphi_n((g_{|\R},...,g_{|\R})(x_1,...,x_n),(g_{|\R},...,g_{|\R})(y_1,...,y_n))&=&\\
\varphi_n((g(x_1),...,g(x_n)),(g(y_1),...,g(y_n)))&=&\\
(g(x_1)-g(y_1))^2+...+(g(x_n)-g(y_n))^2&=&\\
(g(x_1-y_1))^2+...+(g(x_n-y_n))^2&=&\\
g((x_1-y_1)^2)+...+g((x_n-y_n)^2)&=&\\
g((x_1-y_1)^2+...+(x_n-y_n)^2)&=&g((\sqrt{r})^2)=g(r)=r=(\sqrt{r})^2.
\end{eqnarray*}
\par
\noindent
Analogously, if $g:\R \to \F$ is a field homomorphism
then $(g,...,g):{\R}^n \to {\F}^n$ preserves all distances
$\sqrt{r}$ with rational $r \geq 0$.
\vskip 0.2truecm
\par
\noindent
{\bf Conjecture 1.} Each unit-distance preserving mapping from
${\R}^n$ to ${\F}^n$ ($n \geq 2$) has a form $I \circ (g,...,g)$,
where $g:\R \to \F$ is a field homomorphism and $I:{\F}^n \to {\F}^n$
is an affine mapping with orthogonal linear part.
\vskip 0.2truecm
\par
\noindent
{\bf Theorem 1} (\cite{Tyszka200?}). If $x,y \in {\R}^n$ and $|x-y|^2$ is not a rational
number, then there exists $f:{\R}^n \to {\C}^n$ that does not preserve the
distance between $x$ and~$y$ although $f$ preserves all distances
$\sqrt{r}$ with rational $r \geq 0$.
\vskip 0.2truecm
\par
\noindent
{\it Proof.} There exists an isometry $I:{\R}^n \to {\R}^n$
such that $I(x)=(0,0,...,0)$ and $I(y)=(|x-y|,0,...,0)$.
By Proposition 1 there exists
a field automorphism $g:\C \to \C$ such that
$g(|x-y|^2) \neq |x-y|^2$. Thus $g(|x-y|) \neq |x-y|$ and
$g(|x-y|) \neq -|x-y|$. Therefore
$(g_{|\R},...,g_{|\R}): {\R}^n \to {\C}^n$
does not preserve the distance between
$(0,0,...,0) \in {\R}^n$ and $(|x-y|,0,...,0) \in {\R}^n$
although $(g_{|\R},...,g_{|\R})$ preserves all
distances $\sqrt{r}$ with rational $r \geq 0$.
Hence $f:=(g_{|\R},...,g_{|\R}) \circ I: {\R}^n \to {\C}^n$
does not preserve the distance between $x$ and $y$ although
$f$ preserves all distances $\sqrt{r}$ with rational $ r \geq 0$.
\vskip 0.2truecm
\par
Let $A_{n}(\F)$ denote the set of all positive numbers $d$
such that any map $f:{\R}^n \to {\F}^n$ that preserves
unit distance preserves also distance $d$.
The classical Beckman-Quarles theorem states that each unit-distance
preserving mapping from ${\R}^n$ to ${\R}^n$ ($n \geq 2$) is an isometry,
see \cite{Beckman-Quarles}-\cite{Benz1994} and \cite{Everling}.
It means that for each $n \geq 2$ $A_n(\R)=(0,\infty)$.
By Theorem~1 $A_n(\C) \subseteq \left\{d>0: d^2 \in \Q \right\}$.
Let $D_{n}(\F)$ denote the set of all positive numbers~$d$
with the following property:
\begin{description}
\item{($\ast$)}
if $x,y \in {\R}^n$ and $|x-y|=d$ then there exists a finite set $S_{xy}$
with $\left\{x,y \right\} \subseteq S_{xy} \subseteq {\R}^n$ such that any map
$f:S_{xy}\rightarrow {\F}^n$ that preserves unit distance
preserves also the distance between $x$ and $y$.
\end{description}
\par
\noindent
Obviously, $\{1\} \subseteq D_n(\F) \subseteq A_n(\F)$.
From \cite{Tyszka2000} and \cite{Tyszka2001} follows that for each
$n \geq 2$ $D_n(\R)$ is equal to the set of positive
algebraic numbers. The author proved in \cite{Tyszka200?}
that for each $n \geq 2$ $D_n(\C)$ is a dense subset of $(0,\infty)$.
The proof remains valid for $D_n(\F)$.
\vskip 0.2truecm
\par
\noindent
{\bf Theorem 2.} If $n \geq 2$ and $f:{\R}^n \to {\F}^n$
preserves unit distance, then $f$ is injective.
\vskip 0.2truecm
\par
\noindent
{\it Proof.} We know that $D_n(\F)$ is a dense subset of $(0,\infty)$.
Therefore, for any $x,y \in {\R}^n$, $x \neq y$ there exists $z \in {\R}^n$
such that $|z-x| \neq |z-y|$ and $|z-x|,|z-y| \in D_n(\F)$.
All distances in $D_n(\F)$ are preserved by $f$.
Suppose $f(x)=f(y)$. This would imply that
$|z-x|^2=\varphi_n(z,x)= \varphi_n(f(z),f(x))=\varphi_n(f(z),f(y))=\varphi_n(z,y)=|z-y|^2$,
which is a contradiction.
\vskip 0.2truecm
\par
We shall prove that $\{d>0: d^2 \in \Q\} \subseteq D_2(\F)$.
We need the following technical Propositions 2-5.
\vskip 0.2truecm
\par
\noindent
{\bf Proposition~2} (cf. \cite{Blumenthal}, \cite{Borsuk}).
The points $c_{1}=(z_{1,1},...,z_{1,n}),...,
c_{n+1}=(z_{n+1,1},...,z_{n+1,n}) \in {\F}^n$ are affinely dependent
if and only if their Cayley-Menger determinant
$\Delta(c_1,...,c_{n+1}):=$
\footnotesize
$$
\det \left[
\begin{array}{ccccc}
 0  &  1                       &  1                       & ... & 1                       \\
 1  &  0                       & \varphi_n(c_{1},c_{2})   & ... & \varphi_n(c_{1},c_{n+1})\\
 1  & \varphi_n(c_{2},c_{1})   &  0                       & ... & \varphi_n(c_{2},c_{n+1})\\
... & ...                      & ...  	                  & ... & ...                     \\
 1  & \varphi_n(c_{n+1},c_{1}) & \varphi_n(c_{n+1},c_{2}) & ... & 0                       \\
\end{array}\;\right]
$$
\normalsize
\par
\noindent
equals $0$.
\vskip 0.2truecm
\par
\noindent
{\it Proof.} It follows from the equality
$$
\left(
\det \left[
\begin{array}{ccccc}
z_{1,1}   & z_{1,2}   & ... &  z_{1,n}  & 1  \\
z_{2,1}   & z_{2,2}   & ... &  z_{2,n}  & 1  \\
  ...     &  ...      & ... &  ...      & ...\\
z_{n+1,1} & z_{n+1,2} & ... & z_{n+1,n} & 1  \\
\end{array}
\right] \right)^2=
\frac{(-1)^{n+1}}{2^{n}} \cdot \Delta(c_1,...,c_{n+1}).
$$
\vskip 0.2truecm
\par
\noindent
{\bf Proposition~3} (cf. \cite{Blumenthal}, \cite{Borsuk}).
For any points $c_{1},...,c_{n+k} \in {\F}^n$ ($k=2,3,4,...$) their
Cayley-Menger determinant equals $0$ i.e. $\Delta(c_1,...,c_{n+k})=0$.
\vskip 0.2truecm
\par
\noindent
{\it Proof.} Assume that
$c_{1}=(z_{1,1},...,z_{1,n}),...,c_{n+k}=(z_{n+k,1},...,z_{n+k,n})$.
The points $\wt{c}_{1}=(z_{1,1},...,z_{1,n},0,...,0)$,
$\wt{c}_{2}=(z_{2,1},...,z_{2,n},0,...,0)$,...,
$\wt{c}_{n+k}=(z_{n+k,1},...,z_{n+k,n},0,...,0) \in
{\F}^{n+k-1}$ are affinely dependent.
Since $\varphi_n(c_i,c_j)=\varphi_{n+k-1}(\wt{c}_{i},\wt{c}_{j})$
$(1 \leq i \leq j \leq n+k)$ the Cayley-Menger determinant
of points $c_{1},...,c_{n+k}$ is equal to the Cayley-Menger determinant
of points $\wt{c}_{1},...,\wt{c}_{n+k}$ which equals $0$ according to
Proposition~2.
\vskip 0.2truecm
\par
\noindent
{\bf Proposition~4a.}
If $d \in \F$, $d \neq 0$,
$c_{1},...,c_{n+1} \in {\F}^n$ and
$\varphi_n(c_i,c_j)=d^2$ ($1 \leq i<j \leq n+1$),
then the points $c_{1}$,...,$c_{n+1}$ are affinely independent.
\vskip 0.2truecm
\par
\noindent
{\it Proof.} It follows from Proposition 2 because
the Cayley-Menger determinant $\Delta(c_1,...,c_{n+1})=$
$$
\det \left[
\begin{array}{cccccc}
  0  &  1  &  1  & ... &  1  &  1  \\
  1  &  0  & d^2 & ... & d^2 & d^2 \\
  1  & d^2 &  0  & ... & d^2 & d^2 \\ 
 ... & ... & ... & ... & ... & ... \\
  1  & d^2 & d^2 & ... &  0  & d^2 \\
  1  & d^2 & d^2 & ... & d^2 &  0  \\
\end{array}
\right]=(-1)^{n+1}(n+1)d^{2n} \neq 0.
$$
\vskip 0.2truecm
\par
\noindent
{\bf Proposition~4b.} If $a,b \in \F$, $a+b \neq 0$, $z,x,\wt{x} \in {\F}^2$
and $\varphi_2(z,x)=a^2$, $\varphi_2(x,\wt{x})=b^2$,
$\varphi_2(z,\wt{x})=(a+b)^2$, then
$\overrightarrow{zx}=\frac{\textstyle a}{{\textstyle a+b}}\overrightarrow{z\wt{x}}$.
\vskip 0.2truecm
\par
\noindent
{\it Proof.} Since $\varphi_2(z, \wt{x})=(a+b)^2 \neq 0$ we conclude that
$z \neq \wt{x}$.
The Cayley-Menger determinant $\Delta(z,x,\wt{x})=$
$$
\det \left[
\begin{array}{cccc}
0 &   1     &  1   &   1     \\
1 &   0     & a^2  & (a+b)^2 \\
1 &  a^2    &  0   &   b^2   \\
1 & (a+b)^2 & b^2  &   0     \\
\end{array}
\right]
=0,
$$
\par
\noindent
so by Proposition 2 the points $z,x,\wt{x}$ are affinely
dependent. Therefore, there exists $c \in \F$ such that
$\overrightarrow{zx}=c \cdot \overrightarrow{z\wt{x}}$.
Hence $a^2=\varphi_2(z,x)=c^2 \cdot\varphi_2(z,\wt{x})=c^2 \cdot (a+b)^2$.
Thus $c=\frac{\textstyle a}{\textstyle a+b}$
or $c=-\frac{\textstyle a}{\textstyle a+b}$ .
If $c=-\frac{\textstyle a}{\textstyle a+b}$ then
$\overrightarrow{x\wt{x}}=\overrightarrow{z\wt{x}}-\overrightarrow{zx}
=\overrightarrow{z\wt{x}}+\frac{\textstyle a}{\textstyle a+b}\overrightarrow{z\wt{x}}=
\frac{\textstyle 2a+b}{\textstyle a+b}\overrightarrow{z\wt{x}}$.
Hence $b^2=\varphi_2(x,\wt{x})=
\left( \frac{\textstyle 2a+b}{\textstyle a+b} \right)^2 \cdot \varphi_2(z,\wt{x})=
\left( \frac{\textstyle 2a+b}{\textstyle a+b} \right)^2 \cdot (a+b)^2=(2a+b)^2$.
Therefore $0=(2a+b)^2-b^2=4a(a+b)$. Since $a+b \neq 0$ we conclude that
$a=0$, so $c=-\frac{\textstyle a}{\textstyle a+b}=
\frac{\textstyle a}{\textstyle a+b}$ and the proof is completed.
\vskip 0.2truecm
\par
\noindent
{\bf Proposition~5} (cf. [6, Lemma, page 127] for ${\R}^n$).
If $x$, $y$, $c_{0}$, $c_{1}$, $c_{2}$ $\in$ ${\F}^2$,
$\varphi_2(x,c_{0})=\varphi_2(y,c_{0})$,
$\varphi_2(x,c_{1})=\varphi_2(y,c_{1})$,
$\varphi_2(x,c_{2})=\varphi_2(y,c_{2})$
and $c_0$, $c_1$, $c_2$ are affinely independent,
then $x=y$.
\vskip 0.2truecm
\par
\noindent
{\it Proof.} Computing we obtain that the vector
$\overrightarrow{xy}$ is perpendicular
to each of the linearly independent vectors
$\overrightarrow{c_{0}c_{1}}$, $\overrightarrow{c_{0}c_{2}}$.
Thus the vector $\overrightarrow{xy}$ is perpendicular
to every linear combination of vectors
$\overrightarrow{c_{0}c_{1}}$ and $\overrightarrow{c_{0}c_{2}}$.
In particular, the vector
$\overrightarrow{xy}$ is perpendicular
to each of the vectors $[1,0]$, $[0,1]$.
Therefore $\overrightarrow{xy}=0$ and the proof is completed.
\vskip 0.2truecm
\par
\noindent
{\bf Lemma~1}. If $d \in D_2(\F)$ then $\sqrt{3} \cdot d \in D_2(\F)$.
\vskip 0.2truecm
\par
\noindent
{\it Proof.}
Let $ d \in D_2(\F)$, $x,y \in {\R}^2$ and $|x-y|=\sqrt{3} \cdot d$.
Using the notation of Figure~1 we show that
$$
S_{xy}:=S_{y\wt{y}} \cup \bigcup_{i=1}^2 S_{xp_i} \cup \bigcup_{i=1}^2 S_{yp_i} \cup S_{p_1p_2}
\cup
\bigcup_{i=1}^2 S_{x\wt{p_i}}
\cup
\bigcup_{i=1}^2 S_{\wt{y}\wt{p_i}}
\cup S_{\wt{p_1}\wt{p_2}}
$$
satisfies condition $(\ast)$.
\vskip 0.2truecm
\par
\centerline{
\beginpicture
\setcoordinatesystem units <0.5cm, 0.5cm>
\setplotsymbol({.})
\plot 6.00 0.00 0.00 0.00 /
\plot 0.00 0.00 3.00 5.20 /
\plot 3.00 5.20 6.00 0.00 /
\plot 6.00 0.00 9.00 5.20 /
\plot 9.00 5.20 3.00 5.20 /
\plot 9.38 -0.80 4.38 -4.10 /
\plot 4.38 -4.10 4.00 1,88 /
\plot 9.0 5.20 9.38 -0.80 /
\plot 9.00 5.20 4.00 1.88 /
\plot 0.00 0.00 4.38 -4.10 /
\plot 9.38 -0.80 4.00 1.88 /
\setdots
\plot 0.00 0.00 9.0 5.2 /
\plot 9.0 5.20 4.38 -4.10 /
\put {$\wt{p_1}$} at 10 -0.80
\put {$y$} at -0.40 0.00
\put {$p_1$} at 5.80 -0.40
\put {$x$} at 9.00 5.60
\put {$p_2$} at 3.00 5.60
\put {$\wt{y}$} at 4.38 -4.60
\put {$\wt{p_2}$} at 3.60 1.40
\endpicture}
\vskip 0.2truecm
\par
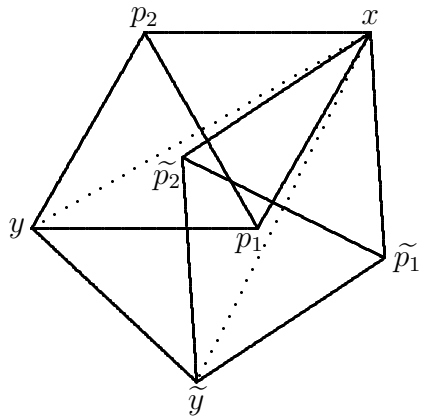
\centerline{Figure~1}
\centerline{$|x-y|=|x-\wt{y}|=\sqrt{3} \cdot d$, $|y-\wt{y}|=|p_1-p_2|=|\wt{p_1}-\wt{p_2}|=d$}
\centerline{$|x-p_i|=|y-p_i|=|x-\wt{p_i}|=|\wt{y}-\wt{p_i}|=d$ ($i=1,2$)}
\vskip 0.2truecm
\par
\noindent
Assume that $f: S_{xy} \rightarrow {\F}^2$
preserves unit distance. Since
$$
S_{xy} \supseteq
S_{y\wt{y}}
\cup
\bigcup_{i=1}^{2}S_{xp_{i}}
\cup\bigcup_ {i=1}^{2}S_{yp_{i}}
\cup
S_{p_{1}p_{2}}
$$
\noindent
we conclude that $f$ preserves the distances between
$y$ and $\wt{y}$,
$x$ and $p_{i}$ ($i=1,2$), $y$ and $p_{i}$ ($i=1,2$),
$p_{1}$ and $p_{2}$.
Hence $\varphi_2(f(y),f(\wt{y}))=
\varphi_2(f(x),f(p_i))=\varphi_2(f(y),f(p_i))=\varphi_2(f(p_1),f(p_2))=d^2$
($i=1,2$).
By Proposition~3 the Cayley-Menger determinant
$\Delta(f(x),f(p_{1}),f(p_{2}),f(y))$ equals $0$ i.e.
\footnotesize
$$
\det\left[
\begin{array}{ccccc}
0 &  1                     &  1                       &  1                       &  1                    \\
1 &  0                     & \varphi_2(f(x),f(p_1))   & \varphi_2(f(x),f(p_2))   & \varphi_2(f(x),f(y))  \\
1 & \varphi_2(f(p_1),f(x)) &  0                       & \varphi_2(f(p_1),f(p_2)) & \varphi_2(f(p_1),f(y))\\
1 & \varphi_2(f(p_2),f(x)) & \varphi_2(f(p_2),f(p_1)) &  0                       & \varphi_2(f(p_2),f(y))\\
1 & \varphi_2(f(y),f(x))   & \varphi_2(f(y),f(p_1))   & \varphi_2(f(y),f(p_2))   &  0                    \\
\end{array}
\right]
=0.
$$
\normalsize
Therefore
$$
\det \left[
\begin{array}{ccccccc}
0 &  1  &  1  &  1  &  1 \\
1 &  0  & d^2 & d^2 &  t \\
1 & d^2 &  0  & d^2 & d^2\\
1 & d^2 & d^2 &  0  & d^2\\
1 &  t  & d^2 & d^2 &  0 \\
\end{array}
\right]
=0
$$
where $t=\varphi_2(f(x),f(y))$. Computing this determinant we obtain
$$
2d^{2}t \cdot (3d^2-t)=0.$$
Therefore
$$t=\varphi_2(f(x),f(y))=\varphi_2(f(y),f(x))=3d^2=|x-y|^2$$
or $$t=\varphi_2(f(x),f(y))=\varphi_2(f(y),f(x))=0.$$
Analogously we may prove that
$$\varphi_2(f(x),f(\wt{y}))=\varphi_2(f(\wt{y}),f(x))=3d^2=|x-\wt{y}|^2$$
or $$\varphi_2(f(x),f(\wt{y}))=\varphi_2(f(\wt{y}),f(x))=0.$$
\newpage
\par
\noindent
If $t=0$ then the points $f(x)$ and $f(y)$ satisfy:
\\
\\
\centerline{$\varphi_2(f(x),f(x))=0=\varphi_2(f(y),f(x))$,}
\\
\centerline{$\varphi_2(f(x),f(p_1))=d^2=\varphi_2(f(y),f(p_1))$,}
\\
\centerline{$\varphi_2(f(x),f(p_2))=d^2=\varphi_2(f(y),f(p_2))$.}
\\
\par
\noindent
By Proposition~4a the points $f(x)$, $f(p_1)$, $f(p_2)$ are affinely
independent. Therefore by Proposition~5 $f(x)=f(y)$ and consequently
$$
d^2=\varphi_2(f(y),f(\wt{y}))=\varphi_2(f(x),f(\wt{y}))
\in \{3d^2,~0 \}.
$$
Since $d^2 \neq 3d^2$ and $d^2 \neq 0$ we conclude that
the case $t=0$ cannot occur. This completes the proof of Lemma~1.
\vskip 0.2truecm
\par
\noindent
{\bf Lemma~2.} If $d \in D_2(\F)$ then $2 \cdot d \in D_2(\F)$.
\vskip 0.2truecm
\par
\noindent
{\it Proof.}
Let $d \in D_2(\F)$, $x,y \in {\R}^2$ and $|x-y|=2 \cdot d$.
Using the notation of Figure~2 we show that
\vskip 0.2truecm
\centerline{$S_{xy}:=
\bigcup \{S_{ab}:a,b\in \{x,y,p_{1},p_{2},p_{3}\},
|a-b|=d \vee |a-b|=\sqrt{3} \cdot d\}$}
\vskip 0.2truecm
\par
\noindent
(where $S_{xp_3}$ and $S_{yp_2}$ are known to exist by Lemma~1)
satisfies condition $(\ast)$.
\vskip 0.2truecm
\par
\noindent
\centerline{
\beginpicture
\normalsize
\setcoordinatesystem units <1mm,1mm>
\setplotsymbol({.})
\plot 0 0 30 0 /
\plot 0 0 -30 0 /
\plot 30 0 15 26 /
\plot 15 26   0 0 /
\plot 0 0 -15 26 /
\plot -15 26 -30 0 /
\plot -15 26 15 26 /
\setdashes
\plot -30 0 15 26 /
\plot 30 0 -15 26 /
\put{$x$} at -30 -2
\put {$p_1$} at 1 -2
\put {$p_2$} at -15 28
\put{$p_3$} at 15 28
\put {$y$} at 30 -2
\endpicture}
\vskip 0.2truecm
\par
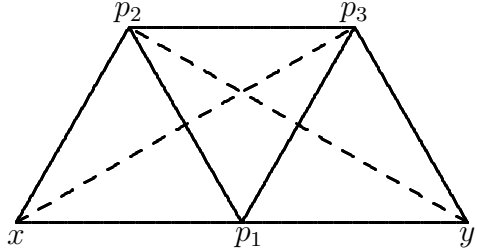
\centerline{Figure~2}
\centerline{$|x-y|=2 \cdot d$}
\centerline{$|p_1-p_2|=|p_1-p_3|=|p_2-p_3|=|x-p_1|=|x-p_2|=|y-p_1|=|y-p_3|=d$}
\centerline{$|x-p_3|=|y-p_2|=\sqrt{3} \cdot d$}
\vskip 0.2truecm
\par
\noindent
Assume that
$f:S_{xy} \rightarrow {\F}^2$ preserves unit distance. Then
$f$ preserves all distances between $p_i$ and $p_j$ $(1 \leq i < j \leq 3)$,
$x$ and $p_i$ $(1 \leq i \leq 3)$, $y$ and $p_i$ $(1 \leq i \leq 3)$.
By Proposition~3 the Cayley-Menger determinant
$\Delta(f(x)$, $f(p_{1})$, $f(p_{2})$, $f(p_{3})$, $f(y))$ equals $0$ i.e.
\par
\footnotesize
$$
\det \left[
\begin{array}{cccccc}
0 &  1                     &  1                       &  1                       &  1                       &  1                    \\
1 &  0                     & \varphi_2(f(x),f(p_1))   & \varphi_2(f(x),f(p_2))   & \varphi_2(f(x),f(p_3))   & \varphi_2(f(x),f(y))  \\
1 & \varphi_2(f(p_1),f(x)) &  0                       & \varphi_2(f(p_1),f(p_2)) & \varphi_2(f(p_1),f(p_3)) & \varphi_2(f(p_1),f(y))\\
1 & \varphi_2(f(p_2),f(x)) & \varphi_2(f(p_2),f(p_1)) &  0                       & \varphi_2(f(p_2),f(p_3)) & \varphi_2(f(p_2),f(y))\\
1 & \varphi_2(f(p_3),f(x)) & \varphi_2(f(p_3),f(p_1)) & \varphi_2(f(p_3),f(p_2)) &  0                       & \varphi_2(f(p_3),f(y))\\
1 & \varphi_2(f(y),f(x))   & \varphi_2(f(y),f(p_1))   & \varphi_2(f(y),f(p_2))   & \varphi_2(f(y),f(p_3))   &  0                    \\
\end{array}
\right]
=0.
$$
\normalsize
Therefore
$$
\det \left[
\begin{array}{cccccc}
0 &  1   &  1  &  1   &  1   &  1  \\
1 &  0   & d^2 & d^2  & 3d^2 &  t  \\
1 & d^2  &  0  & d^2  & d^2  & d^2 \\
1 & d^2  & d^2 &  0   & d^2  & 3d^2\\
1 & 3d^2 & d^2 & d^2  &  0   & d^2 \\
1 &  t   & d^2 & 3d^2 & d^2  &  0  \\
\end{array}
\right]
=0
$$
where $t=\varphi_2(f(x),f(y))$. Computing this determinant we obtain
$$3d^4 \cdot (t-4d^2)^{2}=0.$$
Therefore
$$t=\varphi_2(f(x),f(y))=\varphi_2(f(y),f(x))=4d^2=|x-y|^2.$$
\vskip 0.15truecm
\par
\noindent
{\bf Lemma~3.} If $a,b \in D_2(\F)$ and $a>b$, then $\sqrt{a^2-b^2} \in D_2(\F)$.
\vskip 0.15truecm
\par
\noindent
{\it Proof.}
Let $a, b \in D_2(\F)$, $a>b$, $x,y \in {\R}^2$ and $|x-y|=\sqrt{a^2-b^2}$. Using the notation
of Figure~3 we show that
$$S_{xy}:=S_{xp_1} \cup S_{xp_2} \cup S_{yp_1} \cup S_{yp_2} \cup S_{p_1p_2}$$
\par
\noindent
(where $S_{p_1p_2}$ is known to exist by Lemma~2)
satisfies condition $(\ast)$.
\vskip 0.2truecm
\par
\centerline{
\beginpicture
\setcoordinatesystem units <0.6mm, 0.6mm>
\setplotsymbol({.})
\plot -25 0 0 0 /
\plot 0 0 25 0 /
\plot 25 0 0 42 /
\plot -25 0 0 42 /
\setdots
\plot 0 0 0 42 /
\put {$x$} at 0.5 -4
\put {$y$} at 0.5 45
\put {$p_1$} at -25 -4
\put{$p_2$} at 25 -4
\put {$b$} at -12.5 3.5
\put {$b$} at  12.5 3.5
\put{$a$} at -15 22
\put{$a$} at 16 22
\endpicture}
\vskip 0.2truecm
\par
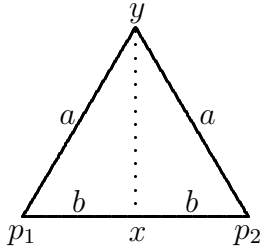
\centerline{Figure~3}
\centerline{$|x-y|=\sqrt{a^2-b^2}$}
\centerline{$|x-p_1|=|x-p_2|=b$, $|y-p_1|=|y-p_2|=a$, $|p_1-p_2|=2b$}
\vskip 0.2truecm
\par
\noindent
Assume that
$f: S_{xy} \rightarrow {\F}^2$ preserves unit distance.
Then $f$ preserves the distances between
$x$ and $p_{i}$ ($i=1,2$), $y$ and $p_{i}$ ($i=1,2$),
$p_1$ and $p_2$.
By Proposition~3 the Cayley-Menger determinant
$\Delta(f(x)$, $f(p_1)$, $f(p_2)$, $f(y))$ equals $0$ i.e.
\footnotesize
$$
\det \left[
\begin{array}{ccccc}
0 &  1                     &  1                       &  1                       &  1                    \\
1 &  0                     & \varphi_2(f(x),f(p_1))   & \varphi_2(f(x),f(p_2))   & \varphi_2(f(x),f(y))  \\
1 & \varphi_2(f(p_1),f(x)) &  0                       & \varphi_2(f(p_1),f(p_2)) & \varphi_2(f(p_1),f(y))\\
1 & \varphi_2(f(p_2),f(x)) & \varphi_2(f(p_2),f(p_1)) &  0                       & \varphi_2(f(p_2),f(y))\\
1 & \varphi_2(f(y),f(x))   & \varphi_2(f(y),f(p_1))   & \varphi_2(f(y),f(p_2))   &  0                    \\
\end{array}
\right]
=0.
$$
\normalsize
Therefore
$$
\det \left[
\begin{array}{ccccc}
0 &  1  &  1   &  1   &  1 \\
1 &  0  & b^2  & b^2  &  t \\
1 & b^2 &  0   & 4b^2 & a^2\\
1 & b^2 & 4b^2 &  0   & a^2\\
1 &  t  & a^2  & a^2  &  0 \\
\end{array}
\right]
=0
$$
where $t=\varphi_2(f(x),f(y))$. Computing this determinant we obtain
$$
-8b^2 \cdot (t+b^2-a^2)^2=0.
$$
Therefore $$t=a^2-b^2=|x-y|^2.$$
\vskip 0.2truecm
\par
\noindent
{\bf Lemma~4.} For each $n \in \{1,2,3,...\}$  $\sqrt{n} \in D_2(\F)$.
\vskip 0.2truecm
\par
\noindent
{\it Proof.} If $n \in \{2,3,4,...\}$ and $\sqrt{n} \in D_2(\F)$,
then $\sqrt{n-1}=\sqrt{(\sqrt{n})^2-1^2} \in D_2(\F)$;
it follows from Lemma 3 and $1 \in D_2(\F)$.
On the other hand, by Lemma 2 all numbers $\sqrt{2^{2k}}=2^k$ ($k=0,1,2,..$) belong
to $D_2(\F)$. These two facts imply that all distances $\sqrt{n}$ ($n=1,2,3,...$) belong
to $D_2(\F)$.
\vskip 0.2truecm
\par
\noindent
From Lemma 4 we obtain:
\vskip 0.2truecm
\par
\noindent
{\bf Lemma~5.} For each $n \in \{1,2,3,...\}$ $n=\sqrt{n^2} \in D_2(\F)$.
\vskip 0.2truecm
\par
\noindent
{\bf Lemma~6.} If $d \in D_2(\F)$ then all distances $\frac{\textstyle d}{\textstyle k}$
($k=2,3,4,...$) belong to $D_2(\F)$.
\vskip 0.2truecm
\par
\noindent
{\it Proof.} Let $d \in D_2(\F)$, $k \in \{2,3,4,...\}$, $x,y \in {\R}^2$ and
$|x-y|=\frac{\textstyle d}{\textstyle k}$. We choose an
integer $e \geq d$. Using the notation of Figure 4 we show that
$$
S_{xy}:=S_{\wt{x}\wt{y}}
\cup
S_{zx}
\cup
S_{x\wt{x}}
\cup
S_{z\wt{x}}
\cup
S_{zy}
\cup
S_{y\wt{y}}
\cup
S_{z\wt{y}}
$$
\par
\noindent
(where sets $S_{zx}$, $S_{x\wt{x}}$, $S_{z\wt{x}}$,
$S_{zy}$, $S_{y\wt{y}}$, $S_{z\wt{y}}$ corresponding
to integer distances are known
to exist by Lemma 5) satisfies condition ($\ast$).
\vskip 0.2truecm
\par
\centerline{
\beginpicture
\setcoordinatesystem units <1.0mm, 1.0mm>
\setplotsymbol({.})
\plot 0 0 -59.3211 9 /
\plot 0 0 -59.3211 -9 /
\plot -59.3211 -9 -59.3211 9 /
\setdots
\plot -19.7737 -3 -19.7737 3 /
\put {$z$} at 2 0
\put {$\wt{y}$} at -59.3211 11.5
\put {$\wt{x}$} at -59.3211 -11.5
\put {$y$} at -19.7737 5
\put {$x$} at -19.7737 -5
\put{$d$} at -61 0
\put {(} at -45 10.70
\put {$k$} at -42 10.25
\put {$-$} at -39 9.80
\put {$1$} at -36 9.35
\put {$)$} at -33 8.90
\put {$e$} at -30 8.45
\put {$e$} at -10 5.45
\put {(} at -45 -10.70
\put {$k$} at -42 -10.25
\put {$-$} at -39 -9.80
\put {$1$} at -36 -9.35
\put {$)$} at -33 -8.90
\put {$e$} at -30 -8.45
\put {$e$} at -10 -5.45
\endpicture}
\vskip 0.2truecm
\par
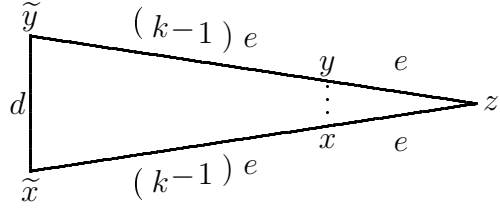
\centerline{Figure 4}
\centerline{$|x-y|=\frac{\textstyle d}{\textstyle k}$}
\vskip 0.2truecm
\par
\noindent
Assume that $f:S_{xy} \to {\F}^2$ preserves unit distance.
Since
$
S_{xy} \supseteq
S_{zx}
\cup
S_{x\wt{x}}
\cup
S_{z\wt{x}}
$
we conclude that:
\par
\noindent
\centerline{$\varphi_2(f(z),f(x))=\varphi_2(z,x)=e^2$,}
\\
\centerline{$\varphi_2(f(x),f(\wt{x}))=\varphi_2(x,\wt{x})=((k-1)e)^2$,}
\\
\centerline{$\varphi_2(f(z),f(\wt{x}))=\varphi_2(z,\wt{x})=(ke)^2$.}
\vskip 0.2truecm
\par
\noindent
By Proposition 4b:
$$
\hspace{5.0cm}
\overrightarrow{f(z)f(x)}=\frac{1}{k}\overrightarrow{f(z)f(\wt{x})}
\hspace{5.0cm}
(1).
$$
Analogously:
$$\hspace{5.0cm}
\overrightarrow{f(z)f(y)}=\frac{1}{k}\overrightarrow{f(z)f(\wt{y})}
\hspace{5.0cm}(2).
$$
By (1) and (2):
$$
\hspace{5cm}
\overrightarrow{f(x)f(y)}=\frac{\textstyle 1}{\textstyle k}\overrightarrow{f(\wt{x})f(\wt{y})}
\hspace{5cm}(3).
$$
Since
$S_{xy} \supseteq S_{\wt{x}\wt{y}}$
we conclude that
$$
\hspace{4.3cm}
\varphi_2(f(\wt{x}),f(\wt{y}))=\varphi_2(\wt{x},\wt{y})=d^2
\hspace{4.3cm}
(4).
$$
By (3) and (4): 
$\varphi_2(f(x),f(y))=\left(\frac{\textstyle d}{\textstyle k}\right)^2$
and the proof is completed.
\vskip 0.2truecm
\par
\noindent
{\bf Theorem 3.} If $x,y \in {\R}^2$ and $|x-y|^2$ is a rational
number, then there exists a finite set $S_{xy}$ with
$\left\{x,y \right\} \subseteq S_{xy} \subseteq {\R}^2$ such that
any map $f:S_{xy}\rightarrow {\F}^2$ that preserves unit distance
preserves also the distance between $x$ and $y$; in other words
\centerline{$\{d>0: d^2 \in \Q \} \subseteq D_2(\F).$
~~~~~~~~~~~~~~~~~~~~~~~~~~~~~~~~~~~~~~~~~~~~~~~~~~~~~~~~~~~~~~~~~~~~~~~~~~~~~}
\vskip 0.2truecm
\par
\noindent
{\it Proof.} We need to prove that
$\sqrt{\frac{\textstyle p}{\textstyle q}} \in D_2(\F)$ for all positive integers $p,q$.
Since $\sqrt{\frac{\textstyle p}{\textstyle q}}=\frac{\sqrt{\textstyle pq}}{\textstyle q}$
the assertion follows from Lemmas 4 and 6.
\vskip 0.2truecm
\par
\noindent
{\bf Conjecture 2.} For each $n \geq 2$
$\{d>0: d^2 \in \Q\} \subseteq D_n(\F)$.
\vskip 0.2truecm
\par
\noindent
As a corollary of Theorem 3 and
$D_2(\C) \subseteq A_2(\C) \subseteq \{d>0: d^2  \in \Q \}$, we get:
\vskip 0.2truecm
\par
\noindent
{\bf Corollary 1.} $D_2(\C)=A_2(\C)=\{d>0: d^2 \in \Q \}$.
\vskip 0.2truecm
\par
The next theorems show which geometric properties are preserved
by unit-distance preserving mappings.
\vskip 0.2truecm
\par
\noindent
{\bf Theorem 4.} If $n \geq 2$, $f:{\R}^n \to {\F}^n$ preserves unit
distance, $x,y \in {\R}^n$ and $x \neq y$, then $\varphi_n(f(x),f(y)) \neq 0$;
in other words $f$ preserves the relation that two points are at non-zero
distance.
\vskip 0.2truecm
\par
\noindent
{\it Proof.} By Theorem 2 $f(x) \neq f(y)$. We know that $D_n(\F)$ is a dense
subset of $(0,\infty)$. In particular, $D_n(\F)$ is unbounded from above.
Therefore, there exist $d \in D_n(\F)$ and $p_1,...,p_{n-1}, p_n \in {\R}^n$
such that $|y-p_i|=d$ ($1 \leq i \leq n-1$) and the points
$x,p_1,...,p_{n-1},p_n$ are at distance $d$ from
each other. Since $f$ preserves distance $d$, we conclude that
$\varphi_n(f(y),f(p_i))=d^2$ ($1 \leq i \leq n-1$) and
$f$ preserves all distances between $x,p_1,...,p_{n-1},p_n$. Thus
$\overrightarrow{f(x)f(p_i)} \cdot \overrightarrow{f(x)f(p_j)}=\frac{d^2}{2}$
($1 \leq i <j \leq n $) and the points $f(x),f(p_1),...,f(p_{n-1}),f(p_n)$
are affinely independent by Proposition 4a. In particular,
the vectors $\overrightarrow{f(x)f(p_i)}$ ($1 \leq i \leq n-1$) are
linearly independent. Assume, on the contrary,
that $\varphi_n(f(x),f(y))=\varphi_n(f(y),f(x))=0$.
It implies that
\vskip 0.2truecm
\par
\noindent
\centerline{$\overrightarrow{f(x)f(y)}$ is perpendicular to each of the
vectors $\overrightarrow{f(x)f(p_i)}$ ($1 \leq i \leq n-1$)~~~~(5)}
\vskip 0.2truecm
\par
\noindent
and the Cayley-Menger determinant
$\Delta(f(x),f(p_1),...,f(p_{n-1}), f(y))=$
$$
\det \left[
\begin{array}{cccccc}
    0   &    1    &       1      & ... &      1       &     1   \\
 {\ONE} & {\ZERO} & {\D}^{\bf 2} & ... & {\D}^{\bf 2} & {\ZERO} \\
    1   &   d^2   &       0      & ... &     d^2      &    d^2  \\
   ...  &   ...   &      ...     & ... &     ...      &    ...  \\
    1   &   d^2   &      d^2     & ... &      0       &    d^2  \\
 {\ONE} & {\ZERO} & {\D}^{\bf 2} & ... & {\D}^{\bf 2} & {\ZERO} \\
\end{array}
\right]
=0,
$$
\par
\noindent
so by Proposition 2 the points $f(x), f(p_1),...,f(p_{n-1}),f(y) \in {\F}^n$
are affinely dependent. Let
$$~~~~~~~~~~~~~~~~~~\overrightarrow{f(x)f(y)}=\alpha_1\overrightarrow{f(x)f(p_1)}+...+\alpha_{n-1}\overrightarrow{f(x)f(p_{n-1})}~~~~~~~~~~~~~~~~~~~~~(6)$$
where $\alpha_1,...,\alpha_{n-1} \in {\F}$ and
\newpage
\begin{eqnarray*}
\overrightarrow{f(x)f(y)}=(t_1,...,t_n) \in {\F}^n,
\\
\overrightarrow{f(x)f(p_1)}=(p_{1,1},...,p_{1,n}) \in {\F}^n,
\\
.~~.~~.~~.~~.~~.~~.~~.~~.~~.~~.~~.~~.~~.~~.~~.~~.~~.
\\
\overrightarrow{f(x)f(p_{n-1})}=(p_{n-1,1},...,p_{n-1,n}) \in {\F}^n.
\end{eqnarray*}
\par
\noindent
By (5):
\begin{eqnarray*}
0=\overrightarrow{f(x)f(p_1)} \cdot \overrightarrow{f(x)f(y)}=p_{1,1} \cdot t_1 + ... + p_{1,n} \cdot t_n,
\\
.~~.~~.~~.~~.~~.~~.~~.~~.~~.~~.~~.~~.~~.~~.~~.~~.~~.~~.~~.~~.~~.~~.~~.~~.~~.
\\
0=\overrightarrow{f(x)f(p_{n-1})} \cdot \overrightarrow{f(x)f(y)}=p_{n-1,1} \cdot t_1 + ... + p_{n-1,n} \cdot t_n.
\end{eqnarray*}
\par
\noindent
By (6):
\begin{eqnarray*}
t_1=\alpha_1 \cdot p_{1,1} + ... + \alpha_{n-1} \cdot p_{n-1,1},
\\
.~~.~~.~~.~~.~~.~~.~~.~~.~~.~~.~~.~~.~~.~~.
\\
t_n=\alpha_1 \cdot p_{1,n} + ... + \alpha_{n-1} \cdot p_{n-1,n}.
\end{eqnarray*}
\par
\noindent
Computing we obtain:
\begin{eqnarray*}
(\overrightarrow{f(x)f(p_1)} \cdot \overrightarrow{f(x)f(p_1)})\alpha_1
+ ... +
(\overrightarrow{f(x)f(p_1)} \cdot \overrightarrow{f(x)f(p_{n-1})})\alpha_{n-1}=0,
\\
.~~.~~.~~.~~.~~.~~.~~.~~.~~.~~.~~.~~.~~.~~.~~.~~.~~.~~.~~.~~.~~.~~.~~.~~.~~.~~.~~.~~.~~.~~.~~.~~.
\\
(\overrightarrow{f(x)f(p_{n-1})} \cdot \overrightarrow{f(x)f(p_1)})\alpha_1
+ ... +
(\overrightarrow{f(x)f(p_{n-1})} \cdot \overrightarrow{f(x)f(p_{n-1})})\alpha_{n-1}=0.
\end{eqnarray*}
\par
\noindent
Therefore:
\begin{eqnarray*}
d^2 \cdot \alpha_1 + \frac{d^2}{2} \cdot \alpha_2 + ... + \frac{d^2}{2} \cdot \alpha_{n-1}=0,
\\
.~~.~~.~~.~~.~~.~~.~~.~~.~~.~~.~~.~~.~~.~~.~~.~~.~~.
\\
\frac{d^2}{2} \cdot \alpha_1 + ... + \frac{d^2}{2} \cdot \alpha_{n-2}+d^2 \cdot \alpha_{n-1}=0.
\end{eqnarray*}
\par
\noindent
Hence $\alpha_1=...=\alpha_{n-1}=0$ and $\overrightarrow{f(x)f(y)}=0$, which
is a contradiction. The proof is completed.
\vskip 0.2truecm
\par
At this point we consider only mappings from ${\R}^2$ to ${\F}^2$.
We omit an easy proof of Proposition 6.
\vskip 0.2truecm
\par
\noindent
{\bf Proposition 6.} If $E,F,C,D \in {\F}^2$, $\varphi_2(E,F) \neq 0$,
$C \neq D$ and
$\varphi_2(E,C)=\varphi_2(F,C)=\varphi_2(E,D)=\varphi_2(F,D)$,
then $\overrightarrow{EC}=\overrightarrow{DF}$ and
$\overrightarrow{FC}=\overrightarrow{DE}$ (see the points
$E,F,C$ and $D$ in Figure 5).
\vskip 0.2truecm
\par
\noindent
{\bf Theorem 5.} If $t \in \Q$, $0<t<1$, $A,B \in {\R}^2$, $A \neq B$,
$C=tA+(1-t)B$ and $f:{\R}^2 \to {\F}^2$ preserves unit distance,
then $f(C)=tf(A)+(1-t)f(B)$.
\vskip 0.2truecm
\par
\noindent
{\it Proof.} We choose $r \in \Q$ such that $r>|AB|$ and
$r<\frac{|AB|}{|1-2t|}$ if $t \neq \frac{1}{2}$. There exists
$D \in {\R}^2$ such that $|AD|=(1-t)r$ and $|BD|=tr$. Let
$E=tA+(1-t)D$ and $F=(1-t)B+tD$, see Figure  5.
\vskip 0.2truecm
\par
\centerline{
\beginpicture
\setcoordinatesystem units <0.8mm, 0.8mm>
\setplotsymbol({.})
\plot -20 0 40 0 /
\plot -6.67 25.82 40 0 /
\plot -23.33 12.9 -20 0 /
\arrow <4mm> [0.1,0.3] from -23.33 12.9 to 0 0
\arrow <4mm> [0.1,0.3] from -6.67 25.82 to 0 0
\arrow <4mm> [0.1,0.3] from -30 38.73 to -23.33 12.9
\arrow <4mm> [0.1,0.3] from -30 38.73 to -6.67 25.82
\put {$C$} at 0 -3
\put {$A$} at -20 -3
\put {$B$} at 40 -3
\put {$D$} at -30 41.73
\put {$F$} at -4.67 28.82
\put {$E$} at -26.33 12.9
\endpicture}
\vskip 0.2truecm
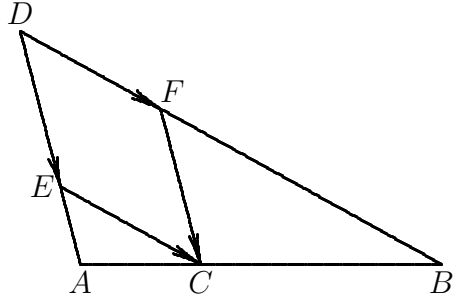
\centerline{Figure 5}
\vskip 0.2truecm
\par
\noindent
The segments $AE$, $ED$, $AD$, $BF$, $FD$, $BD$, $EC$
and $FC$ have rational lengths and
$|EC|=|FC|=|ED|=|FD|=t(1-t)r$.
By Theorem 3 $f$ preserves rational distances. Therefore:
\\
\centerline{$\varphi_2(f(A),f(E))=\varphi_2(A,E)=|AE|^2=((1-t)^2r)^2,$}
\centerline{$\varphi_2(f(E),f(D))=\varphi_2(E,D)=|ED|^2=(t(1-t)r)^2,$}
\centerline{$\varphi_2(f(A),f(D))=\varphi_2(A,D)=|AD|^2=((1-t)r)^2.$}
Since $(1-t)^2r+t(1-t)r=(1-t)r$, by Proposition 4b
$f(E)=tf(A)+(1-t)f(D)$. Analogously
$f(F)=(1-t)f(B)+tf(D)$. Since $C \neq D$, by Theorem 2
$f(C) \neq f(D)$.
Since $E \neq F$, by Theorem 4
$\varphi_2(f(E),f(F)) \neq~0$. By Theorem 3 $f$ preserves
rational distances. Therefore:
\\
\centerline{$\varphi_2(f(E),f(C))=\varphi_2(E,C)=|EC|^2=(t(1-t)r)^2,$}
\centerline{$\varphi_2(f(F),f(C))=\varphi_2(F,C)=|FC|^2=(t(1-t)r)^2,$}
\centerline{$\varphi_2(f(E),f(D))=\varphi_2(E,D)=|ED|^2=(t(1-t)r)^2,$}
\centerline{$\varphi_2(f(F),f(D))=\varphi_2(F,D)=|FD|^2=(t(1-t)r)^2.$}
\vskip 0.2truecm
\par
\noindent
By Proposition 6 $\overrightarrow{f(E)F(C)}=\overrightarrow{f(D)f(F)}$,
so $\overrightarrow{f(A)f(C)}=\overrightarrow{f(A)f(E)}+
\overrightarrow{f(E)f(C)}=\overrightarrow{f(A)f(E)}+\overrightarrow{f(D)f(F)}=
(1-t)\overrightarrow{f(A)f(D)}+(1-t)\overrightarrow{f(D)f(B)}=
(1-t)\overrightarrow{f(A)f(B)}$ and the proof is completed.
\vskip 0.2truecm
\par
It is easy to show that the present form of Theorem 5 implies a
more general form without the assumptions $0<t<1$ and $A \neq B$.
\vskip 0.2truecm
\par
\noindent
{\bf Lemma 7.} If $A,B,C,D \in {\R}^2$,
$\overrightarrow{AB}=\overrightarrow{CD}$,
$|AC|=|BD| \in \Q$ and $f:{\R}^2 \to {\F}^2$ preserves unit
distance, then $\overrightarrow{f(A)f(B)}=\overrightarrow{f(C)f(D)}$.
\vskip 0.2truecm
\par
\noindent
{\it Proof.} There exist $m \in \{0,1,2,...\}$
and $A_0,C_0,A_1,C_1,...,A_m,C_m \in {\R}^2$
such that $A_0=A$, $C_0=C$, $A_m=B$, $C_m=D$ and
for each $i \in \{0,1,...,m-1\}$ $A_iC_iC_{i+1}A_{i+1}$
is a rhombus with a rational side, see Figure 6 where $m=3$.
\vskip 0.2truecm
\par
\centerline{
\beginpicture
\setcoordinatesystem units <1mm, 1mm>
\setplotsymbol({.})
\setdashes
\plot -10 0 0 17.34 /
\plot 0 17.32 14.14 31.46 /
\plot 14.14 31.46 31.46 41.46 /
\plot 10 0 20 17.32 /
\plot 20 17.32 34.14 31.46 /
\plot 34.14 31.46 51.46 41.46 /
\plot -10 0 10 0 /
\plot 0 17.32 20 17.32 /
\plot 14.14 31.46 34.14 31.46 /
\plot 31.46 41.46 51.46 41.46 /
\setsolid
\arrow <4mm> [0.1,0.3] from -10 0 to 31.46 41.46
\arrow <4mm> [0.1,0.3] from 10 0 to 51.46 41.46
\put {$A_0=A$} at -10 -3
\put {$C_0=C$} at 10 -3
\put {$A_1$} at -3 17.32
\put {$C_1$} at 18 20
\put {$A_2$} at 9.14 31.46
\put {$C_2$} at 33 34
\put {$A_3=B$} at 31.46 44.46
\put {$C_3=D$} at 51.46 44.46
\endpicture}
\vskip 0.2truecm
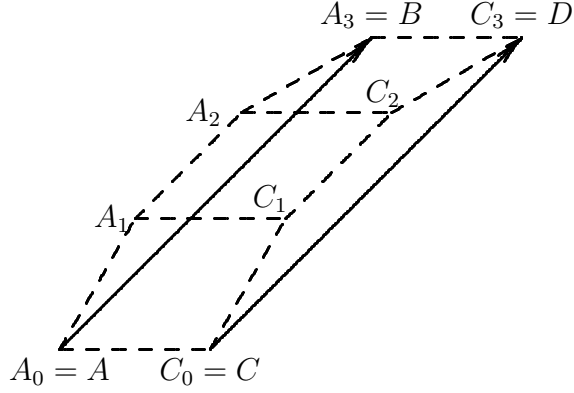
\centerline{Figure 6}
\vskip 0.2truecm
\par
\noindent
By Theorem 3 for each $i \in \{0,1,...,m-1\}$ $f$ preserves
the lenghts of the sides of the rhombus $A_iC_iC_{i+1}A_{i+1}$.
For each $i \in \{0,1,...,m-1\}$ we have:
$\varphi_2(f(A_i),f(C_{i+1})) \neq 0$ (by Theorem 4)
and $f(C_i) \neq f(A_{i+1})$ (by Theorem 2). Therefore,
by Proposition~6 for each $i \in \{0,1,...,m-1\}$
$\overrightarrow{f(A_i)f(A_{i+1})}=\overrightarrow{f(C_i)f(C_{i+1})}$.
Hence $\overrightarrow{f(A)f(B)}=\overrightarrow{f(A_0)f(A_m)}=
\overrightarrow{f(A_0)f(A_1)}+\overrightarrow{f(A_1)f(A_2)}+...+
\overrightarrow{f(A_{m-1})f(A_m)}=\overrightarrow{f(C_0)f(C_1)}+
\overrightarrow{f(C_1)f(C_2)}+...+\overrightarrow{f(C_{m-1})f(C_m)}=
\overrightarrow{f(C)f(D)}$.
\vskip 0.2truecm
\par
\noindent
{\bf Theorem 6.} If $A,B,C,D \in {\R}^2$,
$\overrightarrow{AB}=\overrightarrow{CD}$
and $f:{\R}^2 \to {\F}^2$ preserves unit distance,
then $\overrightarrow{f(A)f(B)}=\overrightarrow{f(C)f(D)}$.
\vskip 0.2truecm
\par
\noindent
{\it Proof.} There exist $E,F \in {\R}^2$ such that
$\overrightarrow{AB}=\overrightarrow{EF}=\overrightarrow{CD}$,
$|AE|=|BF| \in {\Q}$ and $|EC|=|FD| \in {\Q}$, see Figure 7.
\vskip 0.2truecm
\par
\centerline{
\beginpicture
\setcoordinatesystem units <0.7mm, 0.7mm>
\setplotsymbol({.})
\setdashes
\plot 20 0 10 10 /
\plot 10 10 -20 0 /
\plot 30 27 20 37 /
\plot 20 37 -10 27 /
\setsolid
\arrow <4mm> [0.1,0.3] from -20 0 to -10 27
\arrow <4mm> [0.1,0.3] from 20 0 to 30 27
\arrow <4mm> [0.1,0.3] from 10 10 to 20 37
\put {$A$} at -20 -3
\put {$C$} at 20 -3
\put {$E$} at 9 6
\put {$B$} at -13 27
\put {$D$} at 33 27
\put {$F$} at 20 40
\endpicture}
\par
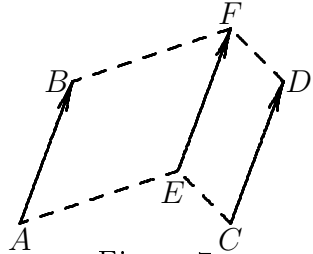
\centerline{Figure 7}
\vskip 0.2truecm
\par
\noindent
By Lemma 7
$\overrightarrow{f(A)f(B)}=\overrightarrow{f(E)f(F)}$
and
$\overrightarrow{f(E)f(F)}=\overrightarrow{f(C)f(D)}$,
so $\overrightarrow{f(A)f(B)}=\overrightarrow{f(C)f(D)}$.
\vskip 0.2truecm
\par
\noindent
As a corollary of Theorems 5 and 6 we get:
\vskip 0.2truecm
\par
\noindent
{\bf Theorem 7.} If $A,B,C,D \in {\R}^2$, $r \in \Q$,
$\overrightarrow{CD}=r\overrightarrow{AB}$ and $f:{\R}^2 \to {\F}^2$
preserves unit distance,
then $\overrightarrow{f(C)f(D)}=r\overrightarrow{f(A)f(B)}$.
\vskip 0.2truecm
\par
\noindent
{\bf Lemma 8.} If $P,Q,X,Y \in {\R}^2$,
$\overrightarrow{PQ}$ is perpendicular to $\overrightarrow{XY}$,
$|XY| \in \Q$ and $f:{\R}^2 \to {\F}^2$ preserves unit distance,
then $\overrightarrow{f(P)f(Q)}$ is perpendicular to
$\overrightarrow{f(X)f(Y)}$.
\vskip 0.2truecm
\par
\noindent
{\it Proof.} There exist $A,B,C,D,E,F \in {\R}^2$ and $r,s \in \Q$
such~that~~$\overrightarrow{PQ}=r\overrightarrow{DE}$,
$\overrightarrow{XY}=s\overrightarrow{AB}$ and the points $A,B,C,D,E,F$
form the configuration from Figure 8; it is a part of Kempe's
linkage for drawing straight lines, see \cite{Rademacher&Toeplitz}.
\vskip 0.2truecm
\par
\centerline{
\beginpicture
\setcoordinatesystem units <0.6mm,0.6mm>
\setplotsymbol({.})
\arrow <4mm> [0.1,0.3] from -40 0 to 40 0
\arrow <4mm> [0.1,0.3] from 5 18 to 5 13.2287
\plot -40 0 5 66.1437 /
\plot 35 39.6832 5 66.1437 /
\plot 40 0 35 39.6832 /
\plot 35 39.6832 5 13.2287 /
\plot 35 39.6832 5 13.2287 /
\plot 5 13.2287 20 0 /
\setdashes
\plot 5 66.1437 5 13.2287 /
\put {$A$} at -40 -4
\put {$B$} at 40 -4
\put {$F$} at 20 -4
\put {$D$} at 3 68.6437
\put {$C$} at 39 39.6862
\put {$E$} at 1.13 13.2287
\endpicture}
\vskip 0.2truecm
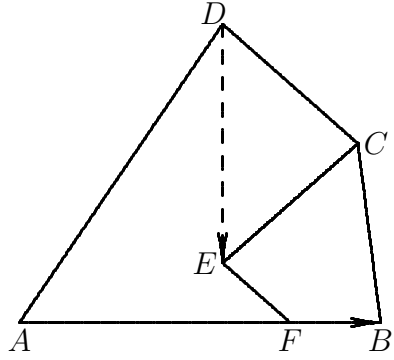
\centerline{Figure 8}
\centerline{$B \neq D$,~~$B \neq E$}
\centerline{$|AB|=|AD|=4$,~~$|CB|=|CD|=|CE|=2$,~~$|AF|=3$,~~$|FB|=|FE|=1$}
\vskip 0.4truecm
\par
\noindent
By Theorem 7 $\overrightarrow{f(P)f(Q)}=r\overrightarrow{f(D)f(E)}$
and $\overrightarrow{f(X)f(Y)}=s\overrightarrow{f(A)f(B)}$, so it
suffices to prove that
$\overrightarrow{f(D)f(E)} \cdot \overrightarrow{f(A)f(B)}=0$.
Let $a=\varphi_2(f(B),f(D))$, $b=\varphi_2(f(A),f(C))$,
$c=\varphi_2(f(B),f(E))$, $d=\varphi_2(f(C),f(F))$,
$e=\varphi_2(f(A),f(E))$.
\vskip 0.2truecm
\par
\noindent
Computing the value of $c$ we obtain:
\par
\noindent
$c=\varphi_2(f(B),f(E))=(\overrightarrow{f(B)f(E)})^2=
(\overrightarrow{f(A)f(E)}-\overrightarrow{f(A)f(B)})^2=$
\par
\noindent
$e-2\overrightarrow{f(A)f(E)} \cdot \overrightarrow{f(A)f(B)}+16$, so
\\
\centerline{~~~~~~~~~~~~~~~~~~~~~~~~~~$\overrightarrow{f(A)f(E)} \cdot \overrightarrow{f(A)f(B)}=8+\frac{1}{2}(e-c)~~~~~~~~~~~~~~~~~~~~~~~~~~~~~~~~$(7).}
\par
\noindent
Computing the value of $a$ we obtain:
\par
\noindent
$a=\varphi_2(f(B),f(D))=(\overrightarrow{f(B)f(D)})^2
=(\overrightarrow{f(A)f(D)}-\overrightarrow{f(A)f(B)})^2=$
\par
\noindent
$16-2\overrightarrow{f(A)f(D)} \cdot \overrightarrow{f(A)f(B)}+16$, so
\par
\noindent
\centerline{~~~~~~~~~~~~~~~~~~~~~~~~~~$\overrightarrow{f(A)f(D)} \cdot \overrightarrow{f(A)f(B)}=
16-\frac{1}{2}a$~~~~~~~~~~~~~~~~~~~~~~~~~~~~~~~~~~~~~(8).}
\newpage
\par
\noindent
The next calculations are based
on Proposition 3 and the observation that distances $1$, $2$, $3$ and $4$
are preserved by $f$.
\vskip 0.2truecm
\par
\noindent
$0=\Delta(f(A),f(B),f(E),f(F))=
\det \left[
\begin{array}{ccccc}
0 &  1 &  1 & 1 & 1 \\
1 &  0 & 16 & e & 9 \\
1 & 16 &  0 & c & 1 \\
1 &  e &  c & 0 & 1 \\
1 &  9 &  1 & 1 & 0 \\
\end{array}
\right]=
-2(e-16+3c)^2$, so
\vskip 0.2truecm
\par
\noindent
\centerline{~~~~~~~~~~~~~~~~~~~~~~~~~~~~~~~~~~~~~~~~~~~$e=16-3c$~~~~~~~~~~~~~~~~~~~~~~~~~~~~~~~~~~~~~~~~~~~~~~~(9).}
\vskip 0.2truecm
\par
\noindent
$0=\Delta(f(A),f(B),f(C),f(F))=
\det \left[
\begin{array}{ccccc}
0 &  1 &  1 &  1 &  1 \\
1 &  0 & 16 &  b &  9 \\
1 & 16 &  0 &  4 &  1 \\
1 &  b &  4 &  0 &  d \\
1 &  9 &  1 &  d &  0 \\
\end{array}
\right]
=-2(b-4d)^2,
$
so $b=4d$.
\vskip 0.4truecm
\par
\noindent
Thus $0=\Delta(f(A),f(B),f(C),f(D))=
\det \left[
\begin{array}{ccccc}
0 &  1 &  1 &  1 &  1 \\
1 &  0 & 16 &  b & 16 \\
1 & 16 &  0 &  4 &  a \\
1 &  b &  4 &  0 &  4 \\
1 & 16 &  a &  4 &  0 \\
\end{array}
\right]=$
\vskip 0.2truecm
\par
\noindent
$\det \left[
\begin{array}{ccccc}
0 &  1 &  1 &  1 &  1 \\
1 &  0 & 16 & 4d & 16 \\
1 & 16 &  0 &  4 &  a \\
1 & 4d &  4 &  0 &  4 \\
1 & 16 &  a &  4 &  0 \\
\end{array}
\right]$
$=-8a(ad+4(d^2-10d+9)).$ By Theorem 4 $a \neq 0$,
so $a=-4\frac{d^2-10d+9}{d}$.
\vskip 0.2truecm
\par
\noindent
$0=\Delta(f(B),f(C),f(E),f(F))=
\det \left[
\begin{array}{ccccc}
0 & 1 & 1 & 1 & 1 \\
1 & 0 & 4 & c & 1 \\
1 & 4 & 0 & 4 & d \\
1 & c & 4 & 0 & 1 \\
1 & 1 & d & 1 & 0 \\
\end{array}
\right]=-2c(cd+d^2-10d+9).$ By Theorem 4 $c \neq 0$, so
$c=-\frac{d^2-10d+9}{d}$. Therefore
\par
\noindent
\centerline{~~~~~~~~~~~~~~~~~~~~~~~~~~~~~~~~~~~~~~~~~~~~~~~~~$a=4c$~~~~~~~~~~~~~~~~~~~~~~~~~~~~~~~~~~~~~~~~~~~~~~~~~~~(10).}
\newpage
\par
\noindent
By (7)-(10):
\vskip 0.2truecm
\par
\noindent
$\overrightarrow{f(D)f(E)} \cdot \overrightarrow{f(A)f(B)}=
(\overrightarrow{f(A)f(E)} - \overrightarrow{f(A)f(D)})
\cdot \overrightarrow{f(A)f(B)}=$
\par
\noindent
$\overrightarrow{f(A)f(E)}\cdot \overrightarrow{f(A)f(B)}-
\overrightarrow{f(A)f(D}) \cdot \overrightarrow{f(A)f(B)}=
8+\frac{1}{2}(e-c)-(16-\frac{1}{2}a)=$
\par
\noindent
$\frac{1}{2}a-\frac{1}{2}c+\frac{1}{2}e-8=
\frac{1}{2}a-\frac{1}{2}c+\frac{1}{2}(16-3c)-8=
\frac{1}{2}a-2c=0$.
\vskip 0.2truecm
\par
\noindent
The proof is completed.
\vskip 0.2truecm
\par
\noindent
{\bf Theorem 8.} If $A,B,C,D \in {\R}^2$, $\overrightarrow{AB}$
and $\overrightarrow{CD}$ are linearly dependent and
$f:{\R}^2 \to {\F}^2$ preserves unit distance,
then $\overrightarrow{f(A)f(B)}$ and 
$\overrightarrow{f(C)f(D)}$ are linearly dependent.
\vskip 0.2truecm
\par
\noindent
{\it Proof.} We choose $X,Y \in {\R}^2$ such that $|XY|=1$ and
both vectors $\overrightarrow{AB}$ and $\overrightarrow{CD}$
are perpendicular to $\overrightarrow{XY}$.
Then $\varphi_2(f(X),f(Y))=1$, so obviously $f(X) \neq f(Y)$.
By Lemma 8 both vectors $\overrightarrow{f(A)f(B)}$
and $\overrightarrow{f(C)f(D)}$ are perpendicular
to $\overrightarrow{f(X)f(Y)}$. These two facts
imply that $\overrightarrow{f(A)f(B)}$ and
$\overrightarrow{f(C)f(D)}$ are linearly dependent.
\vskip 0.2truecm
\par
\noindent
As a corollary of Theorem 8 we get:
\vskip 0.2truecm
\par
\noindent
{\bf Corollary 2.} Unit-distance preserving mappings
from ${\R}^2$ to ${\F}^2$ preserve collinearity of points.
\vskip 0.2truecm
\par
\noindent
{\bf Theorem 9.} If $P,Q,X,Y \in {\R}^2$, $\overrightarrow{PQ}$
is perpendicular to $\overrightarrow{XY}$
and $f:{\R}^2 \to {\F}^2$ preserves unit distance,
then $\overrightarrow{f(P)f(Q)}$ is
perpendicular to $\overrightarrow{f(X)f(Y)}$.
\vskip 0.2truecm
\par
\noindent
{\it Proof.} We may assume that $X \neq Y$. Set
$Z=X+\frac{\overrightarrow{XY}}{|XY|}$.
By Lemma 8 $\overrightarrow{f(P)f(Q)}$ and
$\overrightarrow{f(X)f(Z)}$ are perpendicular.
The points $X$, $Y$ and $Z$ are collinear,
so by Corollary 2 $f(X)$, $f(Y)$ and $f(Z)$ are
collinear. Since $|XZ|=1$, we conclude that
$\varphi_2(f(X),f(Z))=~1$, so $f(X) \neq f(Z)$.
Thus there exists $\alpha \in \F$ such that
$\overrightarrow{f(X)f(Y)}=\alpha\overrightarrow{f(X)f(Z)}$.
Hence $\overrightarrow{f(P)f(Q)} \cdot \overrightarrow{f(X)f(Y)}
=\overrightarrow{f(P)f(Q)} \cdot \alpha \overrightarrow{f(X)f(Z)}=
\alpha\left(\overrightarrow{f(P)f(Q)} \cdot \overrightarrow{f(X)f(Z)}\right)
=\alpha \cdot 0=~0$.
\vskip 0.2truecm
\par
\noindent
{\bf Lemma 9.} If $A,X,Y \in {\R}^2$, $|AX|=|AY|$ and
$f:{\R}^2 \to {\F}^2$ preserves unit distance, then
$\varphi_2(f(A),f(X))=\varphi_2(f(A),f(Y))$.
\vskip 0.2truecm
\par
\noindent
{\it Proof.} Let $Z=\frac{1}{2}X+\frac{1}{2}Y$.
By Theorem 5 $f(Z)=\frac{1}{2}f(X)+\frac{1}{2}f(Y)$.
By Theorem 9 $\overrightarrow{f(A)f(Z)}$ is 
perpendicular to $\overrightarrow{f(Z)f(X)}$
and $\overrightarrow{f(A)f(Z)}$ is perpendicular
to $\overrightarrow{f(Z)f(Y)}$.
Hence
$\varphi_2(f(A),f(X))=
\varphi_2(f(A),f(Z))+\frac{1}{4} \varphi_2(f(X),f(Y))=
\varphi_2(f(A),f(Y))$.
\vskip 0.2truecm
\par
\noindent
{\bf Theorem 10.} If $A,B,C,D \in {\R}^2$, $|AB|=|CD|$ and
$f:{\R}^2 \to {\F}^2$ preserves unit distance,
then $\varphi_2(f(A),f(B))=\varphi_2(f(C),f(D))$.
\vskip 0.2truecm
\par
\noindent
{\it Proof.} We may assume that $|AB|=|CD|>0$.
There exist $m \in \{0,1,2,...\}$ and
points $T_0,T_1,...,T_m \in {\R}^2$ such that
$|AB|=|BT_0|=|T_0T_1|=...=|T_{m-1}T_m|=|T_mC|=|CD|$.
By Lemma 9 $\varphi_2(f(A),f(B))=
\varphi_2(f(B),f(T_0))=
\varphi_2(f(T_0),f(T_1))=~...~=
\varphi_2(f(T_{m-1}),f(T_m))=
\varphi_2(f(T_m),f(C))=
\varphi_2(f(C),f(D))$.
\vskip 0.2truecm
\par
We have already proved the main results that are needed for proving
Conjecture 1 for $n=2$. Such a proof will appear in the next
publication.
\vskip 0.2truecm
\par
{\bf Acknowledgement.} The author wishes to thank Professor Mowaffaq
Hajja for improvement of the proof of Proposition 1.

Apoloniusz Tyszka\\
Technical Faculty\\
Hugo Ko\l{}\l{}\k{a}taj University\\
Balicka 104, 30-149 Krak\'ow, Poland\\
E-mail: {\it rttyszka@cyf-kr.edu.pl}
\end{document}